\documentstyle[twoside, 12pt]{article}

\newtheorem{example}{Example}

\newtheorem{assumption}{Assumption}
\def\qed{ \ \vrule width.2cm height.2cm depth0cm\smallskip}

\newcommand{\ba}{\begin{array}}
\newcommand{\ea}{\end{array}}
\newcommand{\be}{\begin{equation}}
\newcommand{\ee}{\end{equation}}
\newcommand{\bea}{\begin{eqnarray}}
\newcommand{\eea}{\end{eqnarray}}
\newcommand{\beaa}{\begin{eqnarray*}}
\newcommand{\eeaa}{\end{eqnarray*}}

\def\dbR{{\mathop{\rm I\negthinspace R}}}
\def\a{\alpha}
\def\b{\beta}
\def\g{\gamma}
\def\d{\delta}
\def\e{\varepsilon}

\def\l{\lambda}

\def\si{\sigma}
\def\t{\tau}
\def\f{\varphi}
\def\th{\theta}
\def\o{\omega}

\def\G{\Gamma}
\def\D{\Delta}
\def\Th{\Theta}
\def\L{\Lambda}

\def\O{\Omega}
\def\cF{{\cal F}}

\def\no{\noindent}

\def\bs{\bigskip}
\def\q{\quad}
\def\qq{\qquad}

\def\pa{\partial}
\def\cd{\cdot}
\def\cds{\cdots}

\def\td{\nabla}
\def\bF{{\bf F}}
\def\tr{\hbox{\rm tr}}

\def\qed{ \hfill \vrule width.25cm height.25cm depth0cm\smallskip}
\newcommand{\dfnn}{\stackrel{\triangle}{=}}
\newcommand{\basa}{\begin{assumption}}
\newcommand{\easa}{\end{assumption}}

\newcommand{\bas}{\begin{assum}}
\newcommand{\eas}{\end{assum}}

\def\pa{\partial}

 \def\cd{\cdot}
\def\cds{\cdots}

\def\tr{\hbox{\rm tr$\,$}}

\def\dis{\displaystyle}

\def\bF{{\bf F}}

\begin{document}

\newtheorem{thm}{Theorem}[section]
\newtheorem{lem}[thm]{Lemma}
\newtheorem{cor}[thm]{Corollary}
\newtheorem{prop}[thm]{Proposition}
\newtheorem{rem}[thm]{Remark}
\newtheorem{eg}[thm]{Example}
\newtheorem{defn}[thm]{Definition}
\newtheorem{assum}[thm]{Assumption}

\renewcommand {\theequation}{\arabic{section}.\arabic{equation}}
\def\thesection{\arabic{section}}

\title{\bf The Wellposedness of FBSDEs (II)\footnote{This is an old note written in 2005, but was never submitted for publication.}}
\author{{\bf Jianfeng Zhang}\thanks{\noindent
Department of Mathematics, ,University of Southern California,
3620 Vermont Ave, KAP 108, Los Angeles, CA 90089. E-mail:
jianfenz@usc.edu. The author is supported in part by NSF grant
DMS-0403575. }}

\date{}
\maketitle

 \no{\bf Abstract.} This paper is a continuation of \cite{zhang}, in which we established
 the wellposedness result and a comparison theorem for a class of one dimensional
 Forward-Backward SDEs. In this paper we extend the wellposedness result to high dimensional FBSDEs,
 and weaken the key condition in \cite{zhang} significantly. Compared to the existing methods in the literature,
 our result has the following features:
 (i) arbitrary time duration; (ii) random coefficients; (iii) (possibly) degenerate forward diffusion;
 and (iv) no monotonicity condition.

\bs

\no{\bf Keywords}: Forward-backward SDEs, wellposedness.

\bs

\no{\sl MSC 2000 Subject Classifications.} Primary: 60H10


 \section{Introduction and Main Result}
 \setcounter{equation}{0}
 Assume $(\O,\cF,P)$ is a complete probability space, $\cF_0\subset\cF$, and $W$ is a $d$-dimensional
 standard Brownian motion independent of $\cF_0$. Let $\bF\dfnn \{\cF_t\}_{0\le t\le T}$ be the filtration
 generated by $W$ and $\cF_0$, augmented by the null sets as usual. We study the following FBSDE:
 \be
 \label{FBSDE}
 \left\{\ba{lll}
 \dis X_t = X_0 + \int_0^t b(\o,s,\Th_s)ds + \int_0^t \si^*(\o,s,X_s,Y_s) dW_s;\\
 \dis Y_t = g(\o,X_T) +\int_t^T f(\o,s,\Th_s) ds - \int_t^T Z_s dW_s.
 \ea\right.
 \ee
 where $\Th\dfnn (X,Y,Z)$ and $^*$ denotes the transpose. We assume that $X_0\in\cF_0$, $b,\si,f,g$ are
progressively measurable,
  and for any $\th\dfnn(x,y,z)$, $b,\si,f$ are \bF-adapted and $g(\cd,x)\in\cF_T$.
 For simplicity we will always omit the variable $\o$ in
 $b,\si,f,g$.

 The wellposedness of FBSDEs has been studied by many authors (see,
 e.g. \cite{antonelli}, \cite{PT}, \cite{mpy}, \cite{delarue}, \cite{HP},
 \cite{PW}, \cite{yong1}, and \cite{yong2}). We refer the readers to
 \cite{zhang} for a more detailed introduction on the subject.
 Motivated by studying numerical methods for (Markovian) FBSDEs
 (see \cite{CZ}), in \cite{zhang} we proved the following theorem.

 \begin{thm}
 \label{1dmain}
 Assume that all processes are one dimensional; that $b,\si,f,g$ are
 differentiable with respect to $x,y,z$ with uniformly bounded derivatives;
 and that
 \be
 \label{1dkeycondition}
 \pa_y\si\pa_zb = 0;\q \pa_yb + \pa_x\si\pa_zb + \pa_y\si\pa_zf=0.
 \ee
 Denote
 \be
 \label{I0}
 I_0^2\dfnn E\Big\{|X_0|^2 + |g(0)|^2 + \int_0^T [|b(t,0,0,0)|^2 +
 |\si(t,0,0)|^2 + |f(t,0,0,0)|^2]dt\Big\}.
 \ee
 If $I_0^2<\infty$, then FBSDE (\ref{FBSDE}) has a unique solution $\Th$ such
 that
 \be
 \label{norm}
 \|\Th\|^2\dfnn E\Big\{\sup_{0\le t\le T}[|X_t|^2 + |Y_t|^2] + \int_0^T |Z_t|^2 dt\Big\}
 \le CI_0^2.
 \ee
 \end{thm}

 After \cite{zhang} has been accepted for publication, we find that Theorem \ref{1dmain} can be
 improved significantly. In the sequel we assume
 \be
 \label{dimension}
 W\in\dbR^d,\q X,b\in\dbR,\q Y,f,g\in\dbR^n,\q Z\in\dbR^{n\times d},\q \si\in \dbR^d.
 \ee
 Here $W,Y,$ et al are considered as {\it column} vectors. Let $\pa$ denote partial derivatives
 with appropriate dimensions;
 and $|\cd|$ denote the Euclidian norm. For example, $\pa_z b\in \dbR^{n\times d},
 \pa_y\si \in \dbR^{d\times n}$ in an obvious way,
 and $|Y_t|^2=Y_t^*Y_t, |Z_t|^2=\tr(Z_t^*Z_t)$. Our main result is the following theorem.

 \begin{thm}
 \label{main}
 Assume that $b,\si,f,g$ are uniformly Lipschitz continuous in $x,y,z$;
  and that there exists a constant $c>0$ such that
 \be
 \label{keycondition}
  \L^4_t(y) \le -c|\L^3_t(y)|,
 \ee
 for any $y\in\dbR^d$ such that $|y|=1$, where
 \be
 \label{L34y}
 \left.\ba{lll}
 \dis \L^3_t(y) &\dfnn&\dis \sum_{i=1}^n y_i\Big[\tr([\pa_zf^i][\pa_zb]^*)
 -y^*[\pa_zb][\pa_zf^i]^*y  +y^*[\pa_y\si]^*[\pa_zf^i]^*y\Big]\\
 && +[\pa_x\si]^*[\pa_zb]^*y+
 [\pa_yb]y\\
 \dis \L^4_t(y) &\dfnn& |\pa_zb|^2 -|[\pa_zb]^*y|^2 +2y^*[\pa_zb][\pa_y\si]y.
 \ea\right.
 \ee
If $I_0^2<\infty$, then FBSDE (\ref{FBSDE}) has a unique solution $\Th$ such
 that $\|\Th\|^2 \le CI_0^2$, where $C$ depends on $c$ and the Lipschitz constant of the
 coefficients.
 \end{thm}

 We note that we only assume those partial derivatives involved in
 (\ref{L34y}) exist. Moreover, when one part of a product vanishes,
 we do not need to assume the other part to be differentiable. For
 example, if $\pa_zb=0$, then we do not need $\pa_x\si$. In fact, we
 can even weaken (\ref{keycondition}) further by using approximating
 coefficients (see (\ref{approximate}) at below).
 \begin{rem}
 \label{3conditions}
 Following are three sufficient conditions for (\ref{keycondition}):
 \bea
 \label{keycondition1}
 && [\pa_z b][\pa_z b]^* - [\pa_y\si]^*[\pa_zb]^* - [\pa_zb][\pa_y\si] \ge [|\pa_z b|^2 + c]Id_n;\\
 \label{keycondition2}
  && \pa_y b = 0,\q \pa_z b =0,\q [\pa_y\si]^*[\pa_z f^i]^*= 0, i=1,\cds, n;\\
 \label{keycondition3}
 && n=1,\q -[\pa_zb][\pa_y\si]\ge c\Big|\pa_y b + [\pa_zf][\pa_y\si] + [\pa_zb][\pa_x \si]\Big|;
 \eea
 where $Id_n\in \dbR^{n\times n}$ is the $n\times n$ identity
 matrix.
 \end{rem}

  \begin{rem}
 \label{mainrem}
 (i) A necessary condition to ensure (\ref{keycondition1}) is $n\le d$;

 (ii) There are two typical cases for (\ref{keycondition2}). One is that $\pa_y\si=0$, then
 (\ref{FBSDE}) becomes the standard decoupled FBSDE. The other one is that $\pa_z f =0$,
 then (\ref{FBSDE}) becomes
 \be
 \label{FBSDE1}
 \left\{\ba{lll}
 \dis X_t = X_0 + \int_0^t b(s,X_s)ds + \int_0^t \si^*(s,X_s,Y_s) dW_s;\\
 \dis Y_t = g(X_T) +\int_t^T f(s,X_s,Y_s) ds - \int_t^T Z_s dW_s.
 \ea\right.
 \ee
 We note that in this case it is allowed to have $n>d$.
 \end{rem}

 Theorem \ref{main} improves Theorem \ref{1dmain} in three ways. First, there is more freedom on the
 dimensions; second, (\ref{keycondition3}) is obviously much weaker than
 (\ref{1dkeycondition}); and third, we allow the coefficients to be only Lipschitz
 continuous (instead of differentiable). We note that the third feature is not trivial
 because for coefficients satisfying (\ref{keycondition3}) (or (\ref{keycondition1})),
 their molifiers may fail to satisfy so. We would also like to mention that, as in Theorem
 \ref{1dmain}, our result has the following features: 1) $T$ can be arbitrarily large; 2)
 the coefficients are random; 3) $\si$ can be degenerate; 4) no monotonicity condition is
 required.

 However, we should point out that our method does not work when $X$ is high dimensional,
 mainly due to the non-commuting property of matrices
 multiplication. We would leave this
 case for future research.


 \section{Small Time Duration}
 \setcounter{equation}{0}

 In this section we establish some important results for FBSDEs with small time duration
 $T$.
 First we recall a wellposedness result due to Antonelli \cite{antonelli}.
 \begin{lem}
 \label{smallT}
 Assume $b,\si,f$ have a uniform Lipschitz constant $K$,  and $g$ has a uniform Lipschitz constant $K_0$.
 There exist constants $\d_0$ and $C_0$, depending only on $K$ and $K_0$, such that for $T\le\d_0$,
 if $I_0^2<\infty$, then (\ref{FBSDE}) has a unique solution $\Th$ and it holds that
 $
 \|\Th\|\le C_0I_0.
 $
 \end{lem}

  The following lemma, which estimates the $C_0$ at above in terms of $(K,K_0)$,
  is the key step for the proof of Theorem \ref{main}.

 \begin{lem}
 \label{keylem}
  Consider the following linear FBSDE:
  \be
  \label{linearFBSDE}
  \left\{\ba{lll}
 \dis X_t = 1 + \int_0^t B_s ds +  \int_0^t \G_s^* dW_s;\\
 \dis Y_t = GX_T +\int_t^T F_s ds - \int_t^T Z_sdW_s;
 \ea\right.
 \ee
 where
 \beaa
 B_t &=& \a^1_t X_t + \b^1_t Y_t + \tr(\g^1_tZ_t);\\
 \G_t &=& \a^2_t X_t + \b^2_t Y_t;\\
 F_t &=& \a^3_t X_t + \b^3_t Y_t + [\tr(\g^{3,1}_tZ_t),\cds,\tr(\g^{3,n}_tZ_t)]^*;
 \eeaa
 and
 $$
 \a^1_t\in \dbR,\q \b^{1*}_t,\a^3_t\in\dbR^n,\q \a^2_t\in \dbR^d,\q
 \b^2_t,\g^1_t,\g^{3,i}_t\in\dbR^{d\times n},\q \b^3_t\in\dbR^{n\times n}.
 $$
 Assume $|\a^i_t|,|\b^i_t|,|\g^i_t|\le K$, $|G|\le K_0$; and
 \be
 \label{linearkeycondition}
 \L^4_t(y) \le -{1\over K}|\L^3_t(y)|,
 \ee
 for any $y\in\dbR$ such that $|y|=1$,

%

 Let $\d_0$ be as in Lemma \ref{smallT}. There exists a constant $C_K$, depending
 only on $K$ but independent of $K_0$, such that for any $T\le \d_0$, the solution to FBSDE
 (\ref{linearFBSDE}) satisfies
 \be
 \label{K0bar}
 |Y_0|^2\le |\bar K_0|^2 \dfnn [|K_0|^2+1]e^{C_KT}-1.
 \ee
 \end{lem}


 In the sequel we use $C_K$ to denote a generic constant which depends only on $K$ and may vary from line to line.
 Recalling (\ref{L34y}) one can easily check that, for linear FBSDE (\ref{linearFBSDE}), we have
 \be
 \label{linearL34y}
 \left.\ba{lll}
 \dis \L^3_t(y) &=&\dis \sum_{i=1}^n y_i\Big[\tr(\g^{3,i}_t\g^{1*}_t)
 -y^*\g^{1*}_t\g^{3,i}_ty  +y^*\b^{2*}_t\g^{3,i}_ty\Big] +\a^{2*}_t\g^1_ty+
 \b^1_ty\\
 \dis \L^4_t(y) &=& |\g^1_t|^2 -|\g^1_ty|^2 + 2y^*\b^{2*}_t\g^1_ty.
 \ea\right.
 \ee
 We also note that $\tr(AB)=\tr(BA)$ for any matrices $A,B$ with
 appropriate dimensions.

 {\it Proof of Lemma \ref{keylem}.} The proof is quite lengthy, we split it into two steps.

 {\it Step 1.} We first assume $X_t\neq 0$ and formally derive some formulas. Note that
 $$
 B_t\in \dbR;\q \G_t\in\dbR^d;\q F_t\in\dbR^n.
 $$
 Apply Ito's formula, we have
 $$
 dX_t^{-2} = -2X_t^{-3}dX_t + 3X_t^{-4}\G_t^*\G_tdt
 = -2X_t^{-3}\G_t^* dW_t - \Big[2X_t^{-3}B_t -
 3X_t^{-4}\G_t^*\G_t\Big]dt;
 $$
 and
 $$
 d|Y_t|^2 = d(Y_t^*Y_t) = 2Y_t^*dY_t + \tr(Z_tZ_t^*)dt = 2Y_t^*Z_tdW_t -\Big[2Y_t^*F_t -\tr(Z_tZ_t^*)\Big]dt.
 $$
 Denote
 \be
 \label{tildeY}
 \tilde Y_t \dfnn
 Y_tX_t^{-1};\q\tilde\G_t \dfnn \G_t X_t^{-1};\q \tilde Z_t
 \dfnn Z_tX_t^{-1} - \tilde Y_t\tilde\G_t^*;\q d\tilde W_t \dfnn dW_t - [\tilde \G_t+\g^1_t\tilde Y_t]dt.
 \ee
 Recalling that $|Z|^2\dfnn \tr(ZZ^*)$. Then
 \beaa
 && d|\tilde Y_t|^2 = d(|Y_t|^2X_t^{-2}) = X_t^{-2}d|Y_t|^2+|Y_t|^2dX_t^{-2} + d<|Y|^2,X^{-2}>_t\\
 && =  2X_t^{-2}Y_t^*Z_tdW_t -X_t^{-2}[2Y_t^*F_t -|Z_t|^2]dt\\
 &&\q -2|Y_t|^2X_t^{-3}\G_t^* dW_t -|Y_t|^2\Big[2X_t^{-3}B_t - 3X_t^{-4}|\G_t|^2\Big]dt
 -4X_t^{-3}Y_t^*Z_t\G_t dt \\
 &&= \Big[2\tilde Y_t^*Z_tX_t^{-1}-2|\tilde Y_t|^2 \tilde \G_t^*\Big]dW_t\\
 && \q -2\tilde Y_t^*[\a^3_t+\b^3_t\tilde Y_t]dt- 2\sum_{i=1}^n \tilde Y_t^i \tr(\g^{3,i}_tZ_tX_t^{-1})dt
  + |Z_tX_t^{-1}|^2dt\\
 &&\q -2|\tilde Y_t|^2 \Big[\a^1_t + \b^1_t\tilde Y_t + \tr(\g^1_tZ_tX_t^{-1})\Big]dt
 +3|\tilde Y_t|^2|\tilde \G_t|^2dt
   -4\tilde Y_t^*Z_tX_t^{-1}\tilde \G_tdt\\
 &&= 2\tilde Y_t^*\tilde Z_t\Big[d\tilde W_t +  [\tilde \G_t+\g^1_t\tilde Y_t]dt\Big]
 +|\tilde Z_t+ \tilde Y_t\tilde\G_t^*|^2dt\\
 &&\q - 2 \tr\Big([\sum_{i=1}^n \tilde Y_t^i\g^{3,i}_t+|\tilde Y_t|^2 \g^1_t][\tilde Z_t+ \tilde Y_t\tilde\G_t^*]\Big)dt
 -4\tilde Y_t^*[\tilde Z_t+ \tilde Y_t\tilde\G_t^*]\tilde \G_t dt\\
 &&\q - 2\tilde Y_t^*[ \a^3_t + \b^3_t\tilde Y_t]dt
  -2|\tilde Y_t|^2[\a^1_t + \b^1_t\tilde Y_t]dt
    +3|\tilde Y_t|^2|\tilde \G_t|^2dt \\
 &&= 2\tilde Y_t^*\tilde Z_td\tilde W_t + |\tilde Z_t|^2dt
 - 2 \tr\Big(\tilde Z_t[\sum_{i=1}^n \tilde Y_t^i\g^{3,i}_t+|\tilde Y_t|^2\g^1_t - \g^1_t\tilde Y_t\tilde Y^*_t]\Big)dt\\
 &&\q+ |\tilde Y_t\tilde\G_t^*|^2dt
 - 2 \tr\Big([\sum_{i=1}^n \tilde Y_t^i\g^{3,i}_t+|\tilde Y_t|^2 \g^1_t]\tilde Y_t\tilde\G_t^*\Big)dt
  -4\tilde Y_t^*\tilde Y_t\tilde\G_t^*\tilde \G_tdt\\
  &&\q - 2\tilde Y_t^*[ \a^3_t + \b^3_t\tilde Y_t]dt
   -2|\tilde Y_t|^2 [\a^1_t + \b^1_t\tilde  Y_t]dt +3|\tilde Y_t|^2|\tilde \G_t|^2dt \\
 &&\ge 2\tilde Y_t^*\tilde Z_td\tilde W_t - \Big|\sum_{i=1}^n \tilde Y_t^i\g^{3,i}_t
 +|\tilde Y_t|^2\g^1_t - \g^1_t\tilde Y_t\tilde Y^*_t\Big|^2dt\\
 &&\q - 2 \tr\Big([\sum_{i=1}^n \tilde Y_t^i\g^{3,i}_t+|\tilde Y_t|^2 \g^1_t]\tilde Y_t\tilde\G_t^*\Big)dt
   - 2\tilde Y_t^*[ \a^3_t + \b^3_t\tilde Y_t]dt
   -2|\tilde Y_t|^2 [\a^1_t + \b^1_t\tilde  Y_t]dt \\
 &&= 2\tilde Y_t^*\tilde Z_td\tilde W_t
  - \Big[|\sum_{i=1}^n \tilde Y_t^i\g^{3,i}_t|^2 + |\tilde Y_t|^4|\g^1_t|^2 + |\tilde Y_t|^2|\g^1_t\tilde
 Y_t|^2\Big]dt\\
 &&\q +2\Big[\sum_{i=1}^n \tilde Y_t^i\tilde Y^*_t\g^{1*}_t\g^{3,i}_t\tilde Y_t + |\tilde
 Y_t|^2|\g^1_t\tilde Y_t|^2 - |\tilde Y_t|^2\sum_{i=1}^n \tilde
 Y_t^i\tr(\g^{3,i}_t\g^{1*}_t)\Big]dt\\
 &&\q - 2\Big[[\a^2_t+\b^2_t\tilde Y_t]^*[\sum_{i=1}^n \tilde Y_t^i\g^{3,i}_t+|\tilde Y_t|^2 \g^1_t]\tilde Y_t
   +\tilde Y_t^*[ \a^3_t + \b^3_t\tilde Y_t]
   +|\tilde Y_t|^2 [\a^1_t + \b^1_t\tilde  Y_t]\Big]dt \\
  &&\ge 2\tilde Y_t^*\tilde Z_td\tilde W_t - C_K[1+|\tilde Y_t|^2]dt\\
  &&\q - 2\Big[|\tilde Y_t|^2\sum_{i=1}^n \tilde
 Y_t^i\tr(\g^{3,i}_t\g^{1*}_t) -\sum_{i=1}^n \tilde Y_t^i\tilde Y^*_t\g^{1*}_t\g^{3,i}_t\tilde
 Y_t+|\tilde Y_t|^2\a^{2*}_t\g^1_t\tilde Y_t\\
 &&\qq +\sum_{i=1}^n \tilde Y_t^i\tilde Y^*_t\b^{2*}_t\g^{3,i}_t\tilde
 Y_t +|\tilde Y_t|^2\b^1_t\tilde Y_t \Big]dt\\
  &&\q - \Big[|\tilde Y_t|^4|\g^1_t|^2 - |\tilde Y_t|^2|\g^1_t\tilde Y_t|^2 + 2|\tilde
  Y_t|^2 \tilde Y^*_t\b^{2*}_t\g^1_t\tilde Y_t \Big]dt.
\eeaa
 Denote
 $\bar Y_t \dfnn  \tilde Y_t|\tilde Y_t|^{-1}$ when $|\tilde Y_t|\neq 0$, and arbitrary unit vector
 otherwise. Then $|\bar Y_t|=1$ and
 \be
 \label{dYest}
 d|\tilde Y_t|^2 \ge 2\tilde Y^*_t\tilde Z_td\tilde W_t -C_K[1+|\tilde Y_t|^2]dt -\Big[2|\tilde Y_t|^3\L^3_t(\bar Y_t) +|\tilde
 Y_t|^4\L^4_t(\bar Y_t)\Big]dt.
 \ee

\bs

{\it Step 2.} The arguments in this step are similar to those for
Lemma 3.2 in \cite{zhang}, so we will only sketch the main idea.

 Denote
 $$
 \t\dfnn \inf\{t>0: X_t = 0\}\wedge T;\q \t_n\dfnn \inf\{t>0: X_t = {1\over n}\}\wedge T.
 $$
 Then $\t_n\uparrow \t$ and $X_t >0$ for $t\in [0,\t)$. Recall (\ref{tildeY}) for $t\in [0,\t)$.
 By Lemma \ref{smallT} one can easily prove that $|Y_t|\le C_0|X_t|$, and thus
 \be
 \label{tildeYest}
 |\tilde Y_t|\le C_0,\q \forall t\in [0,\t).
 \ee

  By (\ref{linearkeycondition}) we have
 $$
 2|\tilde Y_t|^3\L^3_t(\bar Y_t) + |\tilde Y_t|^4\L^4_t(\bar Y_t) \le -{1\over K}|\L^3_t(\bar Y_t)||\tilde Y_t|^4
 +  2|\tilde Y_t|^3|\L^3_t(\bar Y_t)|
 \le K|\L^3_t(\bar Y_t)||\tilde Y_t|^2.
 $$
 Note that $|\L^3_t(\bar Y_t)|\le C_K$. Thus
 \be
 \label{L34est}
 2|\tilde Y_t|^3\L^3_t(\bar Y_t) + |\tilde Y_t|^4\L^4_t(\bar Y_t) \le C_K|\tilde Y_t|^2.
 \ee
 Then by (\ref{dYest}) one gets
\be
 \label{dYest1}
 d|\tilde Y_t|^2 \ge 2\tilde Y^*_t\tilde Z_td\tilde W_t -C_K[1+|\tilde Y_t|^2]dt;
 \ee
 In light of (\ref{tildeY}) we define
 $$
 M_t = 1 + \int_0^t M_s[\tilde \G_s+\g^1_s\tilde Y_s]^* 1_{\{\t>s\}}
 dW_s;\q L_t= e^{C_Kt},
 $$
 for the $C_K$ in (\ref{dYest1}). By (\ref{tildeYest}) $M$ is a
 martingale. Moreover,
 $$
 d(L_tM_t|\tilde Y_t|^2) \ge (\cds) dW_t - C_KL_tM_tdt,
 $$
 thanks to the obvious fact that $L_t>0, M_t>0$.

 Now for each $n$, we have
 $$
 |\tilde Y_0|^2 \le L_{\t_n}M_{\t_n}|\tilde Y_{\t_n}|^2 - \int_0^{\t_n}(\cds)
 dW_t + C_K\int_0^{\t_n}L_tM_tdt.
 $$
 Thus
 \be
 \label{tildeY0}
 |\tilde Y_0|^2 \le E\Big\{\G_{\t_n}M_{\t_n}|\tilde Y_{\t_n}|^2 + C_K\int_0^{\t_n}L_tM_tdt\Big\}.
 \ee

 On the other hand, if $\t=T$, $|Y_\t|=|Y_T|=|GX_T|=|GX_\t|\le K_0|X_\t|$. If $\t<T$, then $X_\t=0$,
 thus $|Y_\t|\le C_0|X_\t|=0$. Therefore, in both cases it holds that $|Y_\t|\le K_0|X_\t|$. By the same arguments
 as in Lemma 3.2 of \cite{zhang}, one can prove that
 $$
 |\tilde Y_{\t_n}|^2\le |K_0|^2 + C_KE_{\t_n}^{1\over
 2}\{|\t-\t_n|^2\},
 $$
 which, combined with (\ref{tildeY0}), implies that
 \beaa
 |\tilde Y_0|^2 &\le& E\Big\{\G_{\t_n}M_{\t_n}[|K_0|^2 + C_KE_{\t_n}^{1\over
 2}\{|\t-\t_n|^2\}] + C_K\int_0^{\t_n}L_tM_tdt\Big\}\\
 &\le& E\Big\{|K_0|^2\G_{\t_n}M_{\t_n} +
 C_K\int_0^{\t_n}L_tM_tdt\Big\} + C_KE^{1\over
 2}\{|\G_{\t_n}M_{\t_n}|^2\}E^{1\over 2}\{|\t-\t_n|^2\}\\
 &\le& E\Big\{|K_0|^2e^{C_KT}M_{\t_n} +
 C_K\int_0^T e^{C_Kt}M_tdt\Big\} + C_KE^{1\over 2}\{|\t-\t_n|^2\}\\
 &=& |K_0|^2e^{C_KT} + C_K\int_0^T e^{C_Kt}dt + C_KE^{1\over
 2}\{|\t-\t_n|^2\}\\
 &=&|\bar K_0|^2 + C_KE^{1\over 2}\{|\t-\t_n|^2\}.
 \eeaa
 Let $n\to\infty$ and note that $X_0=1$, we prove (\ref{K0bar}).
 \qed

\bs
 We note that estimate (\ref{L34est}) is essential for the wellposedness of
 FBSDEs.

 \begin{example}
  Consider the following one dimensional FBSDE
  \be
  \label{example1}
  \left\{\ba{lll}
  \dis X_t = 1 - \int_0^t Y_s ds;\\
  \dis Y_t = X_T - \int_t^T Z_s dW_s.
  \ea\right.
  \ee
  Then
  $$
  \L^3_t(y) = -y,\q \L^4_t(y)=0.
  $$
  So (\ref{linearkeycondition}) does not hold true. Note that
  $\tilde Y_T = Y_TX_T^{-1} = 1>0$. Actually one can prove in this
  example that $\tilde Y_t>0$ for any $t$, then
  $$
  2|\tilde Y_t|^3\L^3_t(\bar Y_t) + |\tilde Y_t|^4\L^4_t(\bar Y_t) =
  -2\tilde Y_t^3 <0,
  $$
  which implies (\ref{L34est}).
  So we still have $|\tilde Y_0|\le \bar K_0$. Then by using the arguments
  in next section we can show that (\ref{example1}) is wellposeded for arbitrary $T$.
  In fact, (\ref{example1}) satisfies the monotonicity condition in
  \cite{HP}, and thus its wellposedness is already known.
 \end{example}

 We would also like to mention that (\ref{L34est}) is consistent
 with the four step scheme (see \cite{mpy} and \cite{delarue})
 in the following sense. Assume an FBSDE in the four step scheme
 framework has two solutions $\Th^1,\Th^2$. Denote $\tilde Y_t
 =[Y^1_t-Y^2_t][X^1_t-X^2_t]^{-1}$. Note that $Y^i_t=u(t,X^i_t)$ and
 $u$ is uniformly Lipschitz continuous in $x$, where $u$ is the
 solution to the corresponding PDE. Then $\tilde Y_t$ is uniformly bounded and
 thus (\ref{L34est}) holds true.

\bs
 The following result connects FBSDEs (\ref{FBSDE}) and (\ref{linearFBSDE}).
 \begin{cor}
 \label{tdest}
 Assume that all the conditions in Lemma \ref{smallT} as well as (\ref{keycondition}) hold true with $c={1\over K}$.
 Let
 $T\le\d_0$ as in Lemma \ref{smallT}, and
 $\Th^i, i=0,1,$ be the solution to FBSDEs:
 $$
 \left\{\ba{lll}
 \dis X^i_t = x_i + \int_0^t b(s,\Th^i_s)ds + \int_0^t \si^*(s,X^i_s,Y^i_s) dW_s;\\
 \dis Y^i_t = g(X^i_T) +\int_t^T f(s,\Th^i_s) ds - \int_t^T Z^i_s dW_s.
 \ea\right.
 $$
 Then $|Y^1_0-Y^0_0|\le \bar K_0|x_1-x_0|$, where $\bar K_0$ is defined in (\ref{K0bar}).
 \end{cor}
 {\it Proof.}
 We first assume that all the coefficients are differentiable.
 For $0\le\l\le 1$, let $\Th^\l\dfnn (X^\l,Y^\l,Z^\l)$
 and $\td\Th^\l\dfnn (\td X^\l,\td Y^\l,\td Z^\l)$ be the solutions to FBSDEs:
 $$
 \left\{\ba{lll}
 \dis X^\l_t = x_0+\l(x_1-x_0) + \int_0^t b(s,\Th^\l_s)ds + \int_0^t \si^*(s,X^\l_s,Y^\l_s) dW_s;\\
 \dis Y^\l_t = g(X^\l_T) +\int_t^T f(s,\Th^\l_s) ds - \int_t^T Z^\l_s dW_s.
 \ea\right.
 $$
 and
 \be
 \label{Thl}
 \left\{\ba{lll}
 \dis \td X^\l_t = 1 + \int_0^t\Big[\pa_xb(s,\Th^\l_s)\td X^\l_s+ \pa_yb(s,\Th^\l_s)
 \td Y^\l_s+\tr(\pa_zb^*(s,\Th^\l_s)\td Z^\l_s)\Big] ds\\
 \dis\qq\q + \int_0^t [\pa_x\si(s,\Th^\l_s)\td X^\l_s +\pa_y\si(s,\Th^\l_s)\td Y^\l_s]^* dW_s;\\
 \dis \td Y^\l_t = \pa_xg(X^\l_T)\td X^\l_T - \int_t^T\td Z^\l_sdW_s \\
 \dis\qq\q +\int_t^T \Big[\pa_xf(s,\Th^\l_s)\td X^\l_s
  + \pa_yf(s,\Th^\l_s)\td Y^\l_s + \sum_{j=1}^n\tr(\pa_zf^{j*}(s,\Th^\l_s)\td Z^{\l*}_s)\Big] ds ;
 \ea\right.
 \ee
 respectively. One can easily prove that
 $$
 \Th^1_t-\Th^0_t = \int_0^1 {d\over d\l}\Th^\l_t d\l =
 [x_1-x_0]\int_0^1 \td \Th^\l_t d\l.
 $$
 In particular,
 $$
 Y^1_0-Y^0_0 = [x_1-x_0]\int_0^1 \td Y^\l_0 d\l.
 $$
 Note that (\ref{keycondition}) implies (\ref{linearkeycondition})
 for FBSDE (\ref{Thl}). Then by Lemma \ref{keylem} we have $|\td Y^\l_0|\le \bar
 K_0$, and thus
 $$
 |Y^1_0-Y^0_0| \le |x_1-x_0|\int_0^1 |\td Y^\l_0| d\l\le \bar
 K_0|x_1-x_0|.
 $$

 In general case, for any $\e>0$, we may find molifiers
 $(b^\e,\si^\e,f^\e,g^\e)$ such that
 \be
 \label{approximate}
 \L_t^{4,\l,\e}(y) \le -{1\over K}|\L_t^{3,\l,\e}(y)| + \e,
 \ee
 where $\L^{3,\l,\e}$ and $\L^{4,\l,\e}$ are defined in an obvious way, so are other terms
  such as $\Th^{\l,\e}$. Denote $\tilde Y_t^{\l,\e} \dfnn \td
  Y^{\l,\e}_t[\td X^{\l,\e}_t]^{-1}$. By Lemma \ref{smallT} we have
  $|\tilde Y_t^{\l,\e}|\le C_0$ where $C_0$ may depend on $K_0$
  though. Then we have
  \beaa
  &&2|\tilde Y^{\l,\e}_t|^3\L^{3,\l,\e}_t(\bar Y^{\l,\e}_t) + |\tilde Y^{\l,\e}_t|^4\L^{4,\l,\e}_t(\bar
  Y^{\l,\e}_t)\\
  &&  \le 2|\tilde Y^{\l,\e}_t|^3\L^{3,\l,\e}_t(\bar Y^{\l,\e}_t) + |\tilde
  Y^{\l,\e}_t|^4[-{1\over K}|\L_t^{3,\l,\e}(y)| + \e]\\
  &&\le C_K|\tilde Y^{\l,\e}_t|^2 +\e|\tilde Y^{\l,\e}_t|^4
  \le [C_K+\e C_0^2]|\tilde Y^{\l,\e}_t|^2.
  \eeaa
  Now for $\e \le C_0^{-2}$, we know (\ref{L34est}) holds true for
  $\tilde Y^{\l,\e}$, and thus $|Y^{1,\e}_0-Y^{0,\e}_0| \le \bar
  K_0|x_1-x_0|$. Let $\e\to 0$, the lemma follows from the stability result for FBSDEs
  over small time duration (see \cite{antonelli}).
  \qed

 \section{Proof of Theorem \ref{main}}
 \setcounter{equation}{0}

 We now prove Theorem \ref{main} for arbitrarily large $T$. The arguments are exactly the same as in \cite{zhang}.
 So again we will only sketch the main idea. In the sequel we use
 $L_\f$ to denote the smallest Lipschitz constant of a function $\f$.

 {\it Proof.} Let $K$ and $K_0$ be as in Lemma \ref{smallT}.
  By otherwise choosing larger $K$, without loss of generality
  we assume that  $c={1\over K}$ in (\ref{keycondition}). Define
  $\bar K_0$ as in (\ref{K0bar}) (for the arbitrarily large $T$!).
  Let $\d_0$ be a constant as in Lemma  \ref{smallT}, but corresponding to
  $(K,\bar K_0)$ instead of  $(K,K_0)$. Assume $(m-1)\d_0<T\le m\d_0$ for some integer $m$.
  Denote $T_i\dfnn {iT\over m}, i=0,\cds, m$. Define a mapping
 $g_m:\O\times\dbR\to \dbR$ by $g_m(\o,x) \dfnn g(\o,x)$. Now for
 $t\in [T_{m-1}, T_m]$, consider the following FBSDE:
 $$
 \left\{\ba{lll}
 \dis X^m_t = x + \int_{T_{m-1}}^t b(s,\Th^m_s)ds + \int_{T_{m-1}}^t \si^*(s,X^m_s,Y^m_s) dW_s;\\
 \dis Y^m_t = g_m(X^n_{T_m}) +\int_t^{T_m} f(s,\Th^m_s) ds - \int_t^{T_m} Z^m_s dW_s.
 \ea\right.
 $$
 Note that $L_{g_m}\le K_0\le \bar K_0$, by Lemma \ref{smallT} the
 above FBSDE has a unique solution for any $x$. Define
 $g_{m-1}(x)\dfnn Y^m_{T_{m-1}}$. Then for fixed $x$, $g_{m-1}(x)\in
 \cF_{T_{m-1}}$. Moreover, by Corollary \ref{tdest} we have
 $$
 |L_{g_{m-1}}|^2\le |K_1|^2 \dfnn [|K_0|^2+1]e^{C_K(T_m-T_{m-1})}-1 \le |\bar K_0|^2.
 $$

 Next we consider the following FBSDE over $[T_{m-2},T_{m-1}]$:
 $$
 \left\{\ba{lll}
 \dis X^{m-1}_t = x + \int_{T_{m-2}}^t b(s,\Th^{m-1}_s)ds + \int_{T_{m-1}}^t \si^*(s,X^{m-1}_s,Y^{m-1}_s) dW_s;\\
 \dis Y^{m-1}_t = g_{m-1}(X^{m-1}_{T_{m-1}}) +\int_t^{T_{m-1}} f(s,\Th^{m-1}_s) ds - \int_t^{T_{m-1}} Z^{m-1}_s dW_s.
 \ea\right.
 $$
 Similarly we may define $g_{m-2}(x)$ such that
 $$
 |L_{g_{m-2}}|^2\le |K_2|^2 \dfnn [|K_1|^2+1]e^{C_K(T_{m-1}-T_{m-2})}-1 = [|K_0|^2+1]e^{C_K(T_m-T_{m-2})}-1\le
 \bar K_0.
 $$
 Repeat the arguments for $i=m,\cds, 1$, we may define $g_i$ such
 that
 $$
 |L_{g_i}|^2\le |K_{m-i}|^2\dfnn [|K_0|^2+1]e^{C_K(T_m-T_i)}-1\le \bar K_0.
 $$

 Now for any $X_0\in L^2(\cF_0)$, we may construct the solution to FBSDE (\ref{FBSDE})
  piece by piece over subintervals $[T_{i-1},T_i]$ with terminal condition $g_i$, $i=1,\cds, n$.
  Since on each subinterval the solution is unique, we obtain the uniqueness of the solution to FBSDE (\ref{FBSDE}).
  Finally, the estimate $\|\Th\|\le CI_0$ can also be obtained by piece by piece estimates, as done in \cite{zhang}.
  \qed

\bs

Finally we state the stability result whose proof is exactly the
same as in \cite{zhang} and thus is omitted.

 \begin{thm}
 \label{stability1}
 Assume $(b^i,\si^i,f^i,g^i, X^i_0),i=0,1,$ satisfy all the conditions in
 Theorem \ref{main}. Let $\Th^i$ be the corresponding solutions,
 $\D \Th \dfnn \Th^1-\Th^0$, $\D g \dfnn g_1-g_0$, and define other
 terms similarly. Then
 $$
 \|\D\Th\|^2\le CE\Big\{|\D X_0|^2+|\D g(X^1_T)|^2 +
 \int_0^T\Big[|\D b|^2+|\D \si|^2+
 |\D f|^2\Big](t,\Th^1_t)dt\Big\}.
 $$
 \end{thm}

 \begin{cor}
 \label{stability}
 Assume $(b^n,\si^n,f^n,g^n,X^n_0),n=0,1,\cds$ satisfy all the conditions in
 Theorem \ref{main} uniformly; $X^n_0\to X^0_0$ in $L^2$; for $\f=b,\si,f,g$ and for any $(t,\th)$, $\f^n(t,\th)\to
 \f^0(t,\th)$ as $n\to\infty$; and
 $$
 E\Big\{|X^n_0-X_0|^2 + |g^n-g^0|^2(0)+\int_0^T [|b^n-b^0|^2+|\si^n-\si^0|^2+|f^n-f^0|^2](t,0,0,0) dt\Big\}\to 0.
 $$
 Let $\Th^n$ denote the corresponding solutions. Then
 $
 \|\Th^n-\Th^0\|\to 0.
 $
 \end{cor}


\begin{thebibliography}{99}
\bibitem{antonelli}
F. Antonelli, {\it Backward-forward stochastic differential
equations}, {\sl Ann. Appl. Probab.}, 3 (1993), no. 3, 777--793.

\bibitem{CZ}
J. Cvitani\'c and J. Zhang, {\it The steepest decent method for
FBSDEs}, {\sl Electronic Journal of Probability},  10 (2005), 1468-1495.

\bibitem{delarue}
F. Delarue, {\it On the existence and uniqueness of solutions to
FBSDEs in a non-degenerate case}, {\sl Stochastic Process. Appl.},
99 (2002), no. 2, 209--286.

\bibitem{HP}
Y. Hu and S. Peng, {\it Solution of forward-backward stochastic
differential equations}, {\sl Probab. Theory Related Fields}, 103
(1995), no. 2, 273--283.

\bibitem{mpy}
J. Ma, P. Protter, and J. Yong,  {\it Solving forward-backward
stochastic differential equations explicitly - a four step scheme},
{\sl Probab. Theory Relat. Fields.}, 98 (1994), 339-359.

\bibitem{PT}
E. Pardoux and S. Tang, {\it Forward-backward stochastic
differential equations and quasilinear parabolic PDEs}, {\sl
Probab. Theory Related Fields}, 114 (1999), no. 2, 123--150.

\bibitem{PW}
S. Peng and Z. Wu, {\it Fully coupled forward-backward stochastic
differential equations and applications to optimal control}, {\sl
SIAM J. Control Optim.}, 37 (1999), no. 3, 825--843.

\bibitem{yong1}
J. Yong, {\it Finding adapted solutions of forward-backward
stochastic differential equations: method of continuation}, {\sl
Probab. Theory Related Fields}, 107 (1997), no. 4, 537--572.

\bibitem{yong2}
J. Yong, {\it Linear forward-backward stochastic differential
equations}, {\sl Appl. Math. Optim.}, 39 (1999), no. 1, 93--119.

\bibitem{zhang}
J. Zhang, {\it The Wellposedness of FBSDEs}, {\sl Discrete and
Continuous Dynamical Systems -- series B}, 6 (2006), no. 4, 927-940.
\end{thebibliography}
\end{document}